\documentclass[final,3p,times,12pt]{elsarticle}

\usepackage{lineno}%,hyperref
\usepackage{amsmath}
\usepackage{amssymb}
\usepackage{amsthm}
\usepackage{cases}
\usepackage{url}
\usepackage{maplestd2e}
\usepackage[pdfpagelabels,plainpages=false]{hyperref}
\hypersetup{ linkcolor =blue,
 citecolor = cyan,
 urlcolor = blue,
colorlinks = true, bookmarksopen = true, bookmarksnumbered = true,
bookmarksopenlevel = {1}, linktocpage = true,
pdfstartview = FitH,
pdftitle = {A finite sequences of Hahn-type discrete orthogonal polynomials},
pdfsubject
= {article},
pdfauthor = {Jamei et al.},
pdfcreator = {pdftextify}, pdfproducer = {pdftextify}
}

\newtheorem{theorem}{Theorem}

\newtheorem{proposition}[theorem]{Proposition}

\newtheorem{example}[theorem]{Example}

\newtheorem{remark}[theorem]{Remark}
\def \dx {\mathbb{D}_{x}}
\def \sx {\mathbb{S}_{x}}

\newcommand{\hypergeom}[5]{\mbox{$
_#1 F_#2\left( \! \left.
\begin{array}{c}
\multicolumn{1}{c}{\begin{array}{c} #3
\end{array}}\\[1mm]
\multicolumn{1}{c}{\begin{array}{c} #4
            \end{array}}\end{array}
\! \right| \displaystyle{#5}\right) $} }
\newcommand{\qhypergeom}[5]{\mbox{$
_#1 \phi_#2\left(\!  \left.
\begin{array}{c}
\multicolumn{1}{c}{\begin{array}{c} #3
\end{array}}\\[1mm]
\multicolumn{1}{c}{\begin{array}{c} #4
            \end{array}}\end{array}
  \right| \displaystyle{#5}\right) $} }

\journal{ }

\begin{document}

\begin{frontmatter}
 \title{Recurrence equations and their classical orthogonal polynomial solutions on a quadratic or $q$-quadratic lattice }

 \author{Daniel Duviol Tcheutia}
 \ead{duvtcheutia@yahoo.fr}

\address{Institute of Mathematics, University of Kassel, Heinrich-Plett Str. 40, 34132 Kassel, Germany}

\begin{abstract}
Every classical orthogonal polynomial system $p_n(x)$  satisfies a three-term recurrence relation of the type
  \[
    p_{n+1}(x)=(A_nx+B_n)p_n(x)-C_np_{n-1}(x)~ (n=0,1,2,\ldots, p_{-1}\equiv 0),
  \]
  with $C_nA_nA_{n-1}>0$.   Moreover, Favard's theorem  states that the converse is true. A general method to derive the coefficients $A_n$, $B_n$, $C_n$ in terms of the polynomial coefficients of the  divided-difference equations  satisfied by  orthogonal polynomials on a quadratic or $q$-quadratic lattice is recalled. The  Maple implementations  \texttt{rec2ortho} of Koorwinder and Swarttouw or \texttt{retode} of Koepf and Schmersau were developed to identify classical orthogonal polynomials given by their three-term recurrence relation as special functions. The two  implementations \texttt{rec2ortho} and \texttt{retode}  do not handle classical orthogonal polynomials on a quadratic or $q$-quadratic lattice.  In this manuscript, the Maple implementation \texttt{retode} of Koepf and Schmersau is  extended to cover classical orthogonal polynomials on  quadratic or $q$-quadratic lattices and to answer as application an open problem  submitted by  Alhaidari  during the 14th International Symposium on Orthogonal Polynomials, Special Functions and Applications.
\end{abstract}


\end{frontmatter}

\section{Introduction}
Foupouagnigni  showed in \cite{foupouagnigni2008} that     classical orthogonal polynomials
on a quadratic or $q$-quadratic lattice satisfy a second-order divided difference equation of the form
\begin{equation}\label{eq:dxsx}
    \phi(x(s))\dx^2 p_n(x(s))+\psi(x(s))\sx\dx p_n(x(s))+\lambda_n
    p_n(x(s))=0,
\end{equation}
   where $\phi(x)=ax^2+b x+c$, $\psi(x)=dx+e~ (d\neq 0)$, are polynomials of degree at most 2 and of degree one, respectively, the operators $\dx$ and $\sx$ are given by
\[
\dx
f(x(s))=\frac{f(x(s+\frac{1}{2}))-f(x(s-\frac{1}{2}))}{x(s+\frac{1}{2})-x(s-\frac{1}{2})},\quad
\sx f(x(s))= \frac{f(x(s+\frac{1}{2}))+f(x(s-\frac{1}{2}))}{2},
\]
and $x(s)$ is a quadratic or $q$-quadratic lattice defined by \cite{magnus1994}
\begin{equation*}
x(s)=\left\{%
\begin{array}{lll}
    c_1q^s+c_2q^{-s}+c_3 & \textrm{if } &0<q< 1, \\
    c_4s^2+c_5s+c_6 & \textrm{if } &q= 1,
\end{array}%
\right. c_1,\ldots,c_6\in \mathbb{C}.
\end{equation*}
Note that (\ref{eq:dxsx}) is equivalent to a difference or $q$-difference equation of the form
(see \cite[chaps. 9, 14]{KLS})
\begin{equation}\label{eq:sp1sm1}
\lambda_n y(x(s))=B(s)y(x(s+1))-(B(s)+D(s))y(x(s))+D(s)y(x(s-1)),
\end{equation}
with
\[\phi(x(s))=-\frac{1}{2}\Big(x(s+\frac{1}{2})-x(s-\frac{1}{2})\Big)\Big((x(s+1)-x(s))B(s)+(x(s)-x(s-1))D(s)\Big),\]
  \[ \psi(x(s))=(x(s)-x(s+1))B(s)+(x(s)-x(s-1))D(s).\]
Following the work by Foupouagnigni  \cite{foupouagnigni2008}, Njionou Sadjang et al.\ \cite{njionou_et_al_2015}  proved that the Wilson and the continuous dual Hahn polynomials are solutions of a divided-difference equation of the form
\begin{equation}\label{eqWCDH}
    \phi(x)\mathbf{D}_x^2 p_n(x)+\psi(x)\mathbf{S}_x\mathbf{D}_x p_n(x)+\lambda_n
    p_n(x)=0,
\end{equation}
   where the operators $\mathbf{S}_x$ and the Wilson
   operator (see \cite{Cooper2012},    \cite{ismail_stanton2012}) $\mathbf{D}_x $ are  defined by
\[   \mathbf{D}_x f(x)=\frac{f\Big(x+\frac{i}{2}\Big)-f\Big(x-\frac{i}{2}\Big)}{2ix},\quad
   \mathbf{S}_x f(x)=\frac{f\Big(x+\frac{i}{2}\Big)+f\Big(x-\frac{i}{2}\Big)}{2}.
   \]
Using the same approach, Tcheutia et al. \cite{Tcheutia_et_al_2017} derived  a  divided--difference equation of type
 \begin{equation}\label{CHMP}
    \phi(x)\mathbf{\delta}_x^2 y(x)+\psi(x)\mathbf{S}_x\mathbf{\delta}_x y(x)+\lambda_n
    y(x)=0,
\end{equation}
 satisfied by the continuous Hahn and the Meixner--Pollaczek polynomials, where the  difference  operator $\mathbf{\delta}_x $ (see \cite[p. 436]{NIST2010}, compare
\cite[p. 201 and 214]{KLS}, \cite{Njionou2013,njionou_et_al1_2015}, \cite[Equation
(1.15)]{Tratnik_1989}) is defined  as follows:
 \[\mathbf{\delta}_x  f(x)=\frac{f\Big(x+\frac{i}{2}\Big)-f\Big(x-\frac{i}{2}\Big)}{i}.\]
 \eqref{eqWCDH} and \eqref{CHMP} are equivalent to the difference equation (see \cite[chap. 9]{KLS})
 \begin{equation}\label{differnceyxi}
 \lambda_n y(x)=B(x)y(x+i)-(B(x)+D(x))y(x)+D(x)y(x-i),
 \end{equation}
 with
 \[\phi(x)=x((2x+i)B(x)+(2x-i)D(x)),\, \psi(x)=-i((2x+i)B(x)-(2x-i)D(x)),\]
 and
 \[\phi(x)=\frac{1}{2}(B(x)+D(x)),\, \psi(x)=-i(B(x)-D(x)),\]
 respectively.

   The  coefficients of the divided-difference equations given in the forms \eqref{eq:dxsx},  \eqref{eqWCDH} or \eqref{CHMP} can be used for instance  to compute the three--term recurrence relation or some structure formulae, the  inversion coefficients of classical orthogonal polynomials on a quadratic and $q$-quadratic  lattice (see e. g. \cite{FKT2013}, \cite{njionou_et_al_2015}, \cite{Tcheutia_2014}, \cite{Tcheutia_et_al_2017}  and references therein).

 Every classical orthogonal polynomial system $p_n(x)$  satisfies a three-term recurrence relation of the type
  \begin{equation}\label{TTRR}
    p_{n+1}(x)=(A_nx+B_n)p_n(x)-C_np_{n-1}(x) ~(n=0,1,2,\ldots, p_{-1}\equiv 0),
  \end{equation}
  with $C_nA_nA_{n-1}>0$. Moreover, Favard's theorem  states that the converse is true. In Section 2, a general method to derive such three--term recurrence relations for classical orthogonal polynomials on a quadratic or $q$-quadratic lattice in terms of the given polynomials $\phi(x)$ and $\psi(x)$ will be recalled. Alhaidari \cite{Alhaidari} submitted (as open problem during the 14th International Symposium on Orthogonal Polynomials, Special Functions and Applications) two polynomials defined by their three-term recurrence relations and initial values. He was interested by the derivation of their weight functions, generating functions, orthogonality relations, etc.. In order to solve this problem as suggested in the comments by W. Van Assche in \cite{Clarkson_van_Assche}, we use the computer algebra system Maple to identify the polynomials from their recurrence relations, such as in the Maple implementation  \texttt{rec2ortho} of Koorwinder and Swarttouw \cite{Koornwinder_Swarttouw_1998} or \texttt{retode} of Koepf and Schmersau \cite{Koepf_Schmersau_2002}. The two  implementations \texttt{rec2ortho} and \texttt{retode}  do not handle classical orthogonal polynomials on a quadratic or $q$-quadratic lattice.  In Section 3, we extend the Maple implementation \texttt{retode} of Koepf and Schmersau \cite{Koepf_Schmersau_2002} to cover classical orthogonal polynomials on a quadratic or $q$-quadratic lattice and to answer the problem by  Alhaidari \cite{Alhaidari} as application.
\section{Three-term recurrence equation satisfied by classical orthogonal polynomials on a quadratic or $q$-quadratic lattices}\label{sect:TTRR}
 We recall the definition of the \textsl{Pochhammer symbol (or shifted
  factorial)} given by
   \[(a)_m={\Gamma(a+m)\over \Gamma(a)}=a(a+1)(a+2)\cdots (a+m-1),~m=0,1,2,\ldots,
  \]
  and the $q$-Pochhammer symbol
   \[(a;q)_m=\left\{         \begin{array}{ll}
                                                                             \prod\limits_{j=0}^{m-1}(1-aq^j)& \text{ if }\ m=1,2,3,\ldots \\
                                                                             1 &  \text{ if }\  m=0 .                                                                         \end{array}
                                                                         \right.
  \]
 We set the following polynomial basis:
\begin{equation}\label{eq3ch4}
  B_n(a,x)=(aq^s;q)_n(aq^{-s};q)_n=\prod_{k=0}^{n-1}(1-2axq^k+a^2q^{2k}),~n\geq 1,~ B_0(a,x)\equiv 1,
\end{equation}
 where $\displaystyle{ x=x(s)=\cos \theta={q^s+q^{-s}\over 2},~q^s=e^{i\theta}}$;
\begin{equation}\label{eq4ch4}
  \vartheta_n(a,x)=(a+ix)_n(a-ix)_n;
\end{equation}
\begin{equation}\label{eq7ch4}
\left\{
  \begin{array}{ll}
    \xi_n(\gamma,\delta,\mu(x))=(q^{-x};q)_n(\gamma\delta q^{x+1};q)_n=\prod\limits_{k=0}^{n-1}(1+\gamma\delta q^{2k+1}-\mu(x)q^k),~n\geq 1, \\
\xi_0(\gamma,\delta,\mu(x))\equiv 1,
  \end{array}
\right.
\end{equation}
with $\mu(x)=q^{-x}+\gamma\delta q^{x+1}$;
 \begin{equation}  \label{eq8ch4}
{\small \left\{
  \begin{array}{ll}
    \chi_n(\gamma,\delta,\lambda(x))=(-x)_n(x+\gamma+\delta+1)_n=\prod\limits_{k=0}^{n-1}\Big(k(
\gamma+\delta+k+1)-\lambda(x)\Big),~n\geq
1, \\
\chi_0(\gamma,\delta,\lambda(x))\equiv 1,
  \end{array}
\right. }
\end{equation}
 for
$\lambda(x)=x(x+\gamma+\delta+1)$. From the hypergeometric and the basic
hypergeometric representations (see \cite[chaps. 9, 14]{KLS}) of classical orthogonal polynomials
on a quadratic or $q$-quadratic lattice, their  natural bases are
$\{B_n(a,x)\}$, $\{(a+ix)_n\}$, $\{\xi_n(\gamma,\delta,\mu(x))\}$ or
$\{\chi_n(\gamma,\delta,\lambda(x))\}$  whose elements are polynomials of
degree $n$ in the variables $x$, $x$, $\mu(x)$ or $\lambda(x)$,
respectively,  and the basis  $\{\vartheta_n(a,x)\}$ whose elements  are
polynomials of degree $n$ in the variable $x^2$. The operator $\dx$
is appropriate for $B_n(a,x)$,  $\xi_n(\gamma,\delta,\mu(x))$ and
$\chi_n(\gamma,\delta,\lambda(x))$, $\mathbf{\delta}_x$ is appropriate for $\{(a+ix)_n\}$, whereas the corresponding
operator for the basis $\{\vartheta_n(a,x)\}$ is $\mathbf{D}_x$.

Starting from a difference equation of type \eqref{eq:sp1sm1} or \eqref{differnceyxi} given in \cite{KLS}, we deduce  the divided-difference equation of type \eqref{eq:dxsx}, \eqref{eqWCDH} or \eqref{CHMP} satisfied by each classical orthogonal polynomial on a quadratic or $q$-quadratic lattice. Some of them can be found in \cite{foupouagnigni2008}, \cite{njionou_et_al_2015}, \cite{Tcheutia_et_al_2017} and we recall all of them here to make the manuscript self-contained. They will also be recovered using the algorithms implemented in this manuscript.
\subsection{Polynomials expanded in the basis $\{\vartheta_n(\alpha,x)\}$}\label{Wilson}
In this basis  are expanded:\\
1. the Wilson polynomials
\[W_n(x^2;a,b,c,d)=(a+b)_n(a+c)_n(a+d)_n\hypergeom{4}{3}{-n,n+a+b+c+d-1,a+ix,a-ix}
{a+b,a+c,a+d}{1},\]
with
\begin{eqnarray*}
    \phi(x) &=&({x}^{2}-\left( cd+ab+ac+bc+bd+ad \right) {x}+abcd),\\
    \psi(x)&=&(a+b+c+d)x-abc-abd-acd-bcd,~\lambda_n=-n(n-1+a+b+c+d);
  \end{eqnarray*}

2. the continuous dual Hahn polynomials
\[S_n(x^2;a,b,c)=(a+b)_n(a+c)_n\hypergeom{3}{2}{-n,a+ix,a-ix}{a+b,a+c}{1},\]
with
\[\phi(x) =-\left( a+b+c \right) {x}+abc,~\psi(x)={x}-(bc+ab+ac),~\lambda_n=-n.\]
The  procedure to find the coefficients of the recurrence equation \eqref{TTRR} (with $x$ substituted by $ x^2$) in terms of the coefficients $a,~b,~c,~d,~e$ of $\phi(x)$ and $\psi(x)$ is as follows (cf. \cite{FKT2013,  koepfschmersau1998,  Koepf_Schmersau_2002, njionou_et_al_2015,  Tcheutia_et_al_2017}):
\begin{enumerate}
  \item Substitute
  \begin{equation}\label{pn}
  p_n(x):=p_n(x^2)=k_n\vartheta_n(\alpha,x)+k'_n\vartheta_{n}(\alpha,x)+k''_n\vartheta_{n-2}(\alpha,x)+\ldots
  \end{equation}
  in the divided-difference equation \eqref{eqWCDH} (with $x$ substituted by $ x^2$). Next  multiply this equation by  $\vartheta_1(\alpha,x)$ and use the  relations \cite{njionou_et_al_2015}
      \begin{align*}
         &\vartheta_1(\alpha,x)\mathbf{D}^2_x\vartheta_n(\alpha,x)=\eta(n)\eta(n-1)\vartheta_{n-1}(\alpha,x), \\
   & \vartheta_1(\alpha,x)\mathbf{S}_x\mathbf{D}_x \vartheta_n(\alpha,x)=\eta(n)\left(\beta(\alpha+\frac{1}{2},n-1) \vartheta_{n-1}(\alpha,x)+  \vartheta_n(\alpha,x)\right),\\
   &\vartheta_1(\alpha,x)\vartheta_n(\alpha,x)= \nu(\alpha,n ) \vartheta_n(\alpha,x)+  \vartheta_{n+1}(\alpha,x),
      \end{align*}
      with
      \[\eta(n)=n,~\beta(\alpha,n)=-n(n+\alpha-\frac{1}{2}),~\nu(\alpha,n)=-(n^2+2\alpha n). \]
  \item To eliminate the terms $x^2\vartheta_n(\alpha,x)$ and $x^4\vartheta_n(\alpha,x)$,  use twice the relation \cite{njionou_et_al_2015}
  \begin{equation}\label{ttrrxBN1}
   x^2\vartheta_n(\alpha,x)=-(n+\alpha)^2 \vartheta_{n}(\alpha,x)+ \vartheta_{n+1}(\alpha,x).
  \end{equation}
\item\label{step3} Equating the coefficients of $\vartheta_{n+1}(\alpha,x)$ gives $\lambda_n=-n((n-1)a+d)$ in \eqref{eqWCDH}.
Equating the coefficients of $\vartheta_{n}(\alpha,x)$ and $\vartheta_{n-1}(\alpha,x)$ gives $k'_n/k_n,~k''_n/k_n\in \mathbb{Q}(n)$.
\item Substitute the expression of $p_n$ given by \eqref{pn} in the recurrence relation \eqref{TTRR} (with $x$ substituted by $ x^2$) and use \eqref{ttrrxBN1}.  By equating the coefficients of $\vartheta_{n+1}(\alpha,x)$, $\vartheta_{n}(\alpha,x)$, $\vartheta_{n-1}(\alpha,x)$, we get $A_n$, $B_n$ and $C_n$, respectively, in terms of $k_n$, $k'_n$ and $k''_n$.
  \item Substituting the values of $k'_n$ and $k''_n$ given in step \ref{step3} in these equations yields the three unknowns in terms of  $a$, $b$, $c$, $d$, $e$, $n$, $k_{n-1}$, $ k_n$, $ k_{n+1}$.
\end{enumerate}
\begin{proposition}\label{prop1}
Let $p_n(x):=p_n(x^2)=k_n\vartheta_n(\alpha,x)+k'_n\vartheta_{n}(\alpha,x)+k''_n\vartheta_{n-2}(\alpha,x)+\ldots~ (n\in\mathbb{N}_0)$ be a family of polynomial solutions of the divided-difference equation \eqref{eqWCDH}. Then the recurrence equation \eqref{TTRR} (with $x$ substituted by $ x^2$) holds with
\begin{align*}
&\frac{k_n}{k_{n+1}}A_n=1,\\
 & \frac{k_n}{k_{n+1}}B_n=-\frac{n \left( n-1 \right) a \left( 2\,a{n}^{2}-2\,an+4\,nd-2\,b-d \right) -
nd \left(2\,b+d-2\,nd \right) +e \left( 2\,a-d \right)
}{\left(  \left( 2\,n-2 \right) a+d \right)  \left( 2\,an+d \right)},\\
  &\frac{k_{n-1}}{k_{n+1}}C_n=\frac{n \left( an-2\,a+d \right) }{ \left( 2\,an-a+d \right)  \left( 2\,an-3\,a+d \right)  \left( 2\,an-2
\,a+d \right) ^{2}}\times \Big\{\left( n-1 \right) ^{6}{a}^{3}+ \left( n-1 \right) d{b}^{2}\\
&+ \left( -2\, \left( n-1 \right) ^{4}b+
3\, \left( n-1 \right) ^{5}d-4\,c \left( n-1 \right) ^{2} \right) {a}^
{2}+ \left( -2\, \left( n-1 \right) ^{2}{d}^{2}+de \right) b+ \left( -en-c+e \right) {d}^{2}\\
&+ \left(  \left( n-1 \right) ^{2}{b}^{2}-4\, \left( n-1 \right) ^{3
}db+3\, \left( n-1 \right) ^{4}{d}^{2}- \left( n-1 \right)  \left( en+
4\,c-e \right) d-{e}^{2} \right) a+ \left( n-1 \right) ^{3}{d}^{3}
  \Big\}.
\end{align*}
\end{proposition}
\subsection{Polynomials expanded in the basis $\{(\alpha+ix)_n\}$}
The polynomials expanded in this basis are:\\
1. the continuous Hahn polynomials
\[p_n(x;a,b,c,d)=i^n{(a+c)_n(a+d)_n\over n!}\hypergeom{3}{2}{-n,n+a+b+c+d-1,a+ix}{a+c,a+d}{1},\]
with
\[\phi(x)=2\,{x}^{2}-i \left(a+b-c-d \right) x-cd-ab,~\psi(x) =2\left( a+b+c+d \right) x-2\,i \left( ab-cd \right);\]

2. the Meixner-Pollaczek polynomials
\[P_n^{(\lambda)}(x;\theta)={(2\lambda)_n\over n!}e^{in\theta}\hypergeom{2}{1}{-n,\lambda+ix}{2\lambda}{1-e^{-2i\theta}},\]
with
\[\phi(x)= x \cos \left( \theta \right) -\lambda\,\sin \left( \theta
 \right),~\psi(x) =2\left( x\sin \left( \theta \right) +\lambda \cos \left(
\theta \right) \right).
 \]
The action of the operators $\mathbf{\delta}_x $ and $\mathbf{S}_x$ on the  basis $(\alpha+ix)_n$ is given by \cite{Tcheutia_et_al_2017}
\begin{align*}
    (\alpha+ix)\mathbf{\delta}_x ^2(\alpha+ix)_n &= -n(n-1)(\alpha+ix)_{n-1}; \\
  (\alpha+ix)\mathbf{S}_x\mathbf{\delta}_x  (\alpha+ix)_n &= ni(\alpha+ix)_n-\frac{n(n-1)}{2}i(\alpha+ix)_{n-1}; \\
  (\alpha+ix)(\alpha+ix)_n &= (\alpha+ix)_{n+1}-n(\alpha+ix)_n ;\\
  x(\alpha+ix)_n &=-i(\alpha+ix)_{n+1}+i(n+\alpha)(\alpha+ix)_n.
\end{align*}
We suppose that
\[p_n(x)=k_n(\alpha+ix)_n+k'_n (\alpha+ix)_{n-1}+k''_n (\alpha+ix)_{n-2}+\ldots.\]
Using the same approach as in section \ref{Wilson}, it follows that
 $\lambda_n=-n((n-1)a+d)$ in \eqref{CHMP} and the following result holds.
 \begin{proposition}\label{prop:CHMP}
 Let $p_n(x)=k_n(\alpha+ix)_n+k'_n (\alpha+ix)_{n-1}+k''_n (\alpha+ix)_{n-2}+\ldots~ (n\in\mathbb{N}_0)$ be a family of polynomial solutions of the divided-difference equation \eqref{CHMP}. Then the recurrence equation \eqref{TTRR}   is valid with
 \begin{align*}
 &\frac{k_n}{k_{n+1}}A_n=i,\\
 &\frac{k_n}{k_{n+1}}B_n=i\frac {   \left( 2\,b{n}^{2}-2\,bn-2\,e \right) a+d \left( 2\,bn+e \right)}{ \left( 2\,an-2\,a+d \right)  \left( 2\,an+d \right) },\\
 &\frac{k_{n-1}}{k_{n+1}}C_n=n\Big\{-8\,n \left( n-2 \right)  \left( n-1 \right) ^{4}{a}^{5}+ \left( -4\,
 \left( 7\,{n}^{2}-13\,n+2 \right)  \left( n-1 \right) ^{3}d+32\,cn
 \left( n-2 \right)  \left( n-1 \right) ^{2} \right) {a}^{4}\\
 &+ \Big( -8\,n \left( n-2 \right)  \left( n-1 \right) ^{2}{b}^{2}-2\, \left( 19
\,{n}^{2}-34\,n+10 \right)  \left( n-1 \right) ^{2}{d}^{2}+16\,c
 \left( n-1 \right)  \left( 5\,{n}^{2}-9\,n+2 \right) d\\
 &+8\,{e}^{2}n  \left( n-2 \right)  \Big) {a}^{3}+ \Big( -4\, \left( n-1 \right)
 \left( 5\,{n}^{2}-9\,n+2 \right) d{b}^{2}-8\,en \left( n-2 \right) db
+ \left( 72\,{n}^{2}-128\,n+48 \right) {d}^{2}c\\
&- \left( n-1 \right) \left( 5\,n-3 \right)  \left( 5\,n-6 \right) {d}^{3}+ \left( 12\,n-8
 \right) {e}^{2}d \Big) {a}^{2}+ \Big( -4\, \left( 4\,n-3 \right)
 \left( n-1 \right) {d}^{2}{b}^{2}+ \left( -12\,n+8 \right) e{d}^{2}b\\
 &+ \left( 28\,n-24 \right) {d}^{3}c - \left( 8\,n-7 \right)  \left( n-1 \right) {d}^{4}+4\,{e}^{2}{d}^{2} \Big) a+ \left( -4\,n+4 \right) {d}^{3}{b}^{2}-4\,b{d}^{3}e+ \left( 1-n \right) {d}^{5}+4\,c{d}^{4} \Big\} \\
    &\Big/ \Big\{4\left( 2\,an-2\,a+d \right) ^{2} \left( 2\,an-3\,a+d \right)  \left(2\,an+d \right)  \left( 2\,an-a+d \right)\Big\}.
 \end{align*}
 \end{proposition}
\subsection{Polynomials expanded in the basis $\{\ \chi_n(\gamma,\delta,\lambda(x))\}$}\label{sect:racah_DH}
In this basis, we have:\\
1. the Racah polynomials
\[R_n(\lambda(x);\alpha,\beta,\gamma,\delta)=\hypergeom{4}{3}{-n,n+\alpha+\beta+1,-x,
x+\gamma+\delta+1}{\alpha+1,\beta+\delta+1,\gamma+1}{1},~n=0,1,\ldots,N,\]
with
\begin{eqnarray*}
  \phi(\lambda(x)) &=&  \lambda(x)^2+\frac{1}{2} \Big( 2\,\gamma+\gamma
\,\alpha+2\,\gamma\,\delta+\gamma\,\beta+3\,\alpha+2\,\delta+2\,\alpha
\,\beta+\\
&& \delta\,\alpha+3\,\beta+4-\delta\,\beta \Big) \lambda(x) +1/2\, \left( \gamma+1
 \right)  \left( 1+\gamma+\delta \right)  \left( 1+\delta+\beta
 \right)  \left( 1+\alpha \right),\\
  \psi(\lambda(x))&=& \left( 2+\alpha+\beta \right) \lambda(x) + \left( 1+\delta+
\beta \right)  \left( 1+\alpha \right)  \left( \gamma+1 \right);
\end{eqnarray*}
2. the Dual Hahn polynomials
\[R_n(\lambda(x);\gamma,\delta,N)=\hypergeom{3}{2}{-n,-x,x+\gamma+\delta+1}{\gamma+1,-N}{1},~
n=0,1,\ldots,N,\]
with
\begin{eqnarray*}
  \phi(\lambda(x)) = \frac{1}{2} \left( -\gamma-1+2\,N+\delta \right) \lambda(x) +\frac{1}{2}N
 \left( \gamma+1 \right)  \left( 1+\gamma+\delta \right),\,
  \psi(\lambda(x)) = -\lambda(x) +N \left( \gamma+1 \right);
\end{eqnarray*}
We get by direct computations that the action of the operators $\dx$ and $\sx$ on  $\chi_{n}(\gamma,\delta,\lambda(x))$ is given by
\begin{align*}
  &\chi_1(\gamma,\delta,\lambda(x))\ \dx^2\ \chi_n(\gamma,\delta,\lambda(x)) = \eta(n)\eta(n-1)\ \chi_{n-1}(\gamma,\delta,\lambda(x)), \\
  &\chi_{1}(\gamma,\delta,\lambda(x))\ \sx\dx\  \chi_{n}(\gamma,\delta,\lambda(x)) = \eta(n)\left(\beta(\frac{1}{2},\gamma,\delta,n-1)\ \chi_{n-1}(\gamma,\delta,\lambda(x))+\chi_{n}(\gamma,\delta,\lambda(x))\right), \\
  &\chi_{1}(\gamma,\delta,\lambda(x))\ \chi_{n}(\gamma,\delta,\lambda(x))= \nu(\gamma,\delta,n)\ \chi_{n}(\gamma,\delta,\lambda(x))+\chi_{n+1}(\gamma,\delta,\lambda(x)), \\
  &\lambda(x)\ \chi_{n}(\gamma,\delta,\lambda(x))= \mu(\gamma,\delta,n)\ \chi_{n}(\gamma,\delta,\lambda(x))-\chi_{n+1}(\gamma,\delta,\lambda(x)) \\
 \end{align*}
where
\begin{align*}
  &\mu(\gamma,\delta,n)= n(n+\gamma+\delta+1),\  \nu(\gamma,\delta,n) = -n(n+\gamma+\delta+1),\\
  &\beta(a,\gamma,\delta,n) = \frac{-n(2n+\gamma+\delta+2a)}{2}, \ \ \eta(n)=-n.
\end{align*}
We set
\[p_n(x):=p_n(\lambda(x))=k_n\ \chi_n(\gamma,\delta,\lambda(x))+k'_n\ \chi_{n-1}(\gamma,\delta,\lambda(x))+k''_n\  \chi_{n-2}(\gamma,\delta,\lambda(x))+\ldots,\]
use the latter structure relations satisfied by $\chi_n(\gamma,\delta,\lambda(x))$, and proceed as in section \ref{Wilson} to get $\lambda_n=-n((n-1)a+d)$ in  \eqref{eq:dxsx} and the following result.
\begin{proposition}\label{prop_Racah}
 Let $p_n(x):=p_n(\lambda(x))=k_n\ \chi_n(\gamma,\delta,\lambda(x))+k'_n\ \chi_{n-1}(\gamma,\delta,\lambda(x))+k''_n\  \chi_{n-2}(\gamma,\delta,\lambda(x))+\ldots~ (n\in\mathbb{N}_0)$ be a family of polynomial solutions of the divided-difference equation \eqref{eq:dxsx} (where $\phi(x)\leftarrow \phi(\lambda(x))$ and $\psi(x) \leftarrow \psi(\lambda(x))$). Then the recurrence equation \eqref{TTRR} (with $x\leftarrow \lambda(x)$)   is valid with
\begin{align*}
&\frac{k_n}{k_{n+1}}A_n=-1,\\
   & \frac{k_n}{k_{n+1}}B_n=-\frac{an \left( n-1 \right)  \left( 2\,a{n}^{2}-2\,an+4\,dn+2\,b-d \right) +
2\,bdn+n \left( 2\,n-1 \right) {d}^{2}+de-2\,ae
}{\left( 2\,an-2\,a+d \right)  \left( 2\,an+d \right)},\\
   &\frac{k_{n-1}}{k_{n+1}}C_n=-\frac{n \left( an-2\,a+d \right)}{4\, \left( 2\,an-a+d \right)  \left( 2\,an-3\,a+d \right)  \left( 2\,a
n-2\,a+d \right) ^{2}
}\times \Big\{\left( 4\,c+ \left( -4\,n+4 \right) e \right) {d}^{2}\\
&+ \left( -8\, \left( n-1 \right) ^{4}b+
4\, \left( n-1 \right) ^{3} \left( 2\,{\delta}^{2}+4\,\delta\,\gamma+2 \,{\gamma}^{2}-3\,{n}^{2}+4\,\delta+4\,\gamma+6\,n-1 \right) d+16\,c
 \left( n-1 \right) ^{2} \right) {a}^{2}\\
 &+ \Big( -4\, \left( n-1  \right) ^{2}{b}^{2}-16\, \left( n-1 \right) ^{3}db+ \left( n-1  \right) ^{2} \left(  5\,{\delta}^{2}+10\,\delta\,\gamma+5\,{\gamma}^{2}-12\,{n}^{2}+10\,\delta+10\,\gamma+24\,n-7 \right) {d}^{2}\\
 &+ \left( \left( 16\,n-16 \right) c-4\,e \left( n-1 \right) ^{2} \right) d+4\,{e}^{2} \Big) a+4\, \left( n-1 \right) ^{4} \left( n+\delta+\gamma \right)  \left( -n+2+\delta+\gamma \right) {a}^{3} + \left( -4\,n+4 \right) d{b}^{2}\\
&+ \left( -8\, \left(n-1 \right) ^{2}{d}^{2}-4\,de \right) b+ \left( n-1 \right)  \left( -2
\,n+3+\delta+\gamma \right)  \left( 2\,n-1+\delta+\gamma \right) {d}^{3}
    \Big\}.
\end{align*}
\end{proposition}
\subsection{Polynomials expanded in the basis $\{B_n(\alpha,x)\}$}\label{sect_Askey_Wilson}
The following polynomial families are expanded in the basis $\{B_n(\alpha,x)\}$:\\
 1. the Askey-Wilson polynomials defined  by
 \[p_n(x;a,b,c,d|q)={(ab,ac,ad;q)_n\over a^n}\qhypergeom{4}{3}{q^{-n},abcdq^{n-1},ae^{i\theta},ae^{-i\theta}}{ab,ac,ad}{q;q},~x=\cos\theta,\]
  for which
   \begin{eqnarray*}
  \phi(x) &=&\frac{\left( q-1 \right) ^{2}}{4{q}^{3/2}} \Big(2( abcd+1 )   x^2- \left(a+b+c+d+acd+abc+abd+bcd \right) x \\
  && +ab+cd+bd+bc+ac+ad-abcd-1 \Big), \\
  \psi(x) &=& \frac { \left( q-1 \right)}{2q}\Big(  2\left( abcd-1 \right) x +a+b+c+d-abc-abd-acd-bcd ) \Big);
\end{eqnarray*}
  2. the continuous dual $q$-Hahn polnomials
  \[p_n(x;a,b,c|q)={(ab,ac;q)_n\over a^n}\qhypergeom{3}{2}{q^{-n},ae^{i\theta},ae^{-i\theta}}{ab,ac}{q;q},~x=\cos\theta,\]
  with
  \begin{eqnarray*}
  \phi(x) &=& \frac { \left( q-1 \right) ^{2}}{4{q}^{3/2}}\Big( 2 x^2- \left( a+b+c+abc \right) x
  +  bc+ac+ab-1 \Big), \\
  \psi(x) &=& \frac { \left( q-1 \right)}{2q}\Big( -2 x +a+b+c-abc \Big);
\end{eqnarray*}
  3. the continuous $q$-Hahn polynomials
  {\small\[p_n(x;a,b,c,d;q)={(abe^{2i\hat{\theta}},ac,ad;q)_n\over
(ae^{i\hat{\theta}})^n}\qhypergeom{4}{3}{q^{-n},abcdq^{n-1},ae^{i(\theta+2\hat{\theta})}
,ae^{-i\theta}}
{abe^{2i\hat{\theta}},ac,ad}{q;q},~x=\cos(\theta+\hat{\theta}),\]}
with
\begin{eqnarray*}
\phi(x)& =& \frac{\left( q-1 \right) ^{2}}{4q^{3/2}} \Big( 2\, \left( abcd+1 \right)
 x^2-\frac { \left( d+bcd+a{t}^{2
}+abd{t}^{2}+abc{t}^{2}+c+acd+b{t}^{2} \right)  }{t}x  \\
   &&+\frac {ac{t}^{2}+bd{t}^{2}-abcd{t}^{2}+bc{t}^{2}+cd-{t}^{2}+ab{t}^
{4}+ad{t}^{2}}{{t}^{2}} \Big), ~t=e^{i\hat{\theta}}, \\
\psi(x)& =&  \frac { \left( q-1 \right)}{2q}\Big( 2 \left( abcd-1 \right) x +\frac{\left( c+d-acd+b{t}^{2}+a{t}^
{2}-bcd-abc{t}^{2}-abd{t}^{2} \right)}{t} \Big) ;
\end{eqnarray*}
  4. the  Al-Salam-Chihara polynomials
  \[Q_n(x;a,b|q)={(ab;q)_n\over a^n}\qhypergeom{3}{2}{q^{-n},ae^{i\theta},ae^{-i\theta}}{ab,0}{q;q},~x=\cos\theta,\]
  with
  \[\phi(x)=\frac { \left( q-1 \right) ^{2}}{4{q}^{3/2}}\Big(2x^2- \left( b+a \right) x  + ab-1 \Big)
   ,\ \psi(x)= \frac { \left( q-1 \right)}{2q} \Big( -2x  + b+a \Big);\]
  5. the $q$-Meixner-Pollaczek polynomials
  \[P_n(x;a|q)=a^{-n}e^{-in\hat{\theta}}{(a^2;q)_n\over (q;q)_n}\qhypergeom{3}{2}{q^{-n},ae^{i(\theta+2\hat{\theta})},ae^{-i\theta}}
{a^2,0}{q;q},~x=\cos(\theta+\hat{\theta}),\]
with $(t=e^{i\hat{\theta}})$
 \[\phi(x)=\frac{(q-1)^2}{4q^{3\over 2}}\Big(2x^2-\frac{a(1+t^2)}{t}x+a^2-1\Big),\, \psi(x)=\frac{q-1}{2q}\Big(-2x+\frac{a(1+t^2)}{t}\Big);\]
   6. the continuous $q$-Jacobi polynomials
   \[P_n^{(\alpha,\beta)}(x|q)={(q^{\alpha+1};q)_n\over (q;q)_n}\qhypergeom{4}{3}{q^{-n},q^{n+\alpha+\beta+1},q^{{\alpha\over 2}+{1\over 4}}e^{i\theta},q^{{\alpha\over 2}+{1\over 4}}e^{-i\theta}}{q^{\alpha+1},-q^{\alpha+\beta+1\over 2},-q^{\alpha+\beta+2\over 2}}{q;q},~x=\cos\theta,\]
   with ($p^2=q$)
   \begin{eqnarray*}
   &&\phi(x)= \frac{(p^2-1)^2}{4p^3}\Big(2(p^{2\alpha+2\beta+4}+1)x^2+\sqrt{p}(p+1)(p^\alpha-p^\beta)(p^{\alpha+\beta+2}-1)x \\
  &&+p^{2\alpha+2}+p^{2\beta+2}-p^{\alpha+\beta+1}-p^{\alpha+\beta+3}-2p^{\alpha+\beta+2}-p^{2\alpha+2\beta+4}-1 \Big),\\
  &&\psi(x)=\frac{p^2-1}{2p^2}\Big(2(p^{2\alpha+2\beta+4}-1)x+\sqrt{p}(p+1)(p^\alpha-p^\beta)(p^{\alpha+\beta+2}+1) \Big);
\end{eqnarray*}
   7. the  continuous $q$-ultraspherical/Rogers polynomials
   \[C_n(x;\beta|q)={(\beta^2;q)_n\over (q;q)_n}\beta^{-{n\over 2}}\qhypergeom{4}{3}{q^{-n},\beta^2q^n,\beta^{1\over 2}e^{i\theta},\beta^{1\over 2}e^{-i\theta}}{\beta q^{1\over 2},-\beta,-\beta q^{1\over 2}}{q;q},~x=\cos \theta,\]
   with
   \[ \phi(x)=\left( q-1 \right) ^{2} \left(  \left( 2\,{\beta}^{2}q+2 \right) {x}^{2}- \left( \beta+1 \right)  \left( \beta q+1 \right)  \right),~\psi(x)=4\, \left( q-1 \right) x \left( {\beta}^{2}q-1 \right) \sqrt {q};
\]
   8. the continuous big $q$-Hermite polynomials
   \[H_n(x;a|q)=a^{-n}\qhypergeom{3}{2}{q^{-n},ae^{i\theta},ae^{-i\theta}}{0,0}{q;q},~x=\cos\theta,\]
   with
   \begin{eqnarray*}
  \phi(x) &=& \frac{(q-1)^2}{4q^{3\over 2}}\Big(2x^2-ax-1 \Big), \   \psi(x) = \frac{q-1}{2q}\Big(-2x+a \Big);
\end{eqnarray*}
   9. the continuous $q$-Laguerre polynomials
   \[P_n^{(\alpha)}(x|q)={(q^{\alpha+1};q)_n\over (q;q)_n}\qhypergeom{3}{2}{q^{-n},q^{{\alpha\over 2}+{1\over 4}}e^{i\theta},q^{{\alpha\over 2}+{1\over 4}}e^{-i\theta}}{q^{\alpha+1},0}{q;q},~x=\cos\theta,\]
   with
   \begin{eqnarray*}
   &&\phi(x)=\frac{(p^2-1)^2}{4p^3}\Big(2x^2-p^{\alpha+{1\over 2}}(p+1)x +p^{2\alpha+2}-1\Big),\\
   && \psi(x)= \frac{p^2-1}{2p^2}\Big(-2x+p^{\alpha+{1\over 2}}(p+1)\Big).
\end{eqnarray*}
    The  procedure to find the coefficients of the recurrence equation \eqref{TTRR} in terms of the coefficients $a,~b,~c,~d,~e$ of $\phi(x)$ and $\psi(x)$ is as in section \ref{Wilson}:
\begin{enumerate}
  \item Substitute
  \begin{equation}\label{Pn}
  p_n(x)=k_nB_n(\alpha,x)+k'_nB_{n-1}(\alpha,x)+k''_nB_{n-2}(\alpha,x)+\ldots
  \end{equation}
  in the divided-difference equation \eqref{eq:dxsx}. Next  multiply this equation by  $B_1(\alpha,x)$ and use the  relations \cite{FKT2013}
      \begin{align}
         B_1(\alpha, x)\dx^2B_n(\alpha, x) &= \eta(\alpha,n)\eta(\alpha\sqrt{q},n-1)B_{n-1}(\alpha, x), \label{basisbn1} \\
         B_1(\alpha,x)\sx\dx B_n(\alpha,x) &= \eta(\alpha,n)\left(\beta_1(\alpha\sqrt{q},n-1)B_{n-1}(\alpha, x)+\beta_2(n-1)B_n(\alpha,x)\right),\label{basisbn2}\\
         B_1(\alpha,x)B_n(\alpha,x) &= \nu_1(\alpha,n)B_n(\alpha, x)+\nu_2(n)B_{n+1}(\alpha, x),\label{basisbn3}
      \end{align}
      with
      \begin{align*}
         &\eta(\alpha,n)=\frac{2\alpha(1-q^n)}{q-1},\ \ \beta_1(\alpha,n) = \frac{1}{2}(1-\alpha^2q^{2n-1})(1-q^{-n}),  \\
         & \beta_2(n)=\frac{1}{2}+\frac{1}{2q^n},\  \nu_1(\alpha,n) = (1-q^{-n})(1-\alpha^2q^n),\ \ \nu_2(n)=q^{-n}.
      \end{align*}
  \item To eliminate the terms $xB_n(\alpha,x)$ and $x^2B_{n}(\alpha,x)$,   use  the relations \cite{FKT2013}
\begin{align}
   xB_n(\alpha,x) &= \mu_1(\alpha,n)B_n(\alpha, x)+\mu_2(\alpha,n)B_{n+1}(\alpha, x),\label{ttrrxBN}  \\
   x^2B_n(\alpha,x)&=\mu_1^2(\alpha,n)B_n(\alpha,x)+
    \mu_2(\alpha,n)(\mu_1(\alpha,n)+\mu_1(\alpha,n+1))B_{n+1}(\alpha,x)\nonumber\\
&+\mu_2(\alpha,n)\mu_2(\alpha,n+1)B_{n+2}(\alpha,x), \label{ttrrx2BN}
\end{align}
with
\[\mu_1(\alpha,n) = \frac{1+\alpha^2q^{2n}}{2 \alpha q^n},\ \ \mu_2(\alpha,n)=\frac{-1}{2\alpha q^n}.\]
\item Equating
the coefficients of $B_{n+1}(\alpha,x)$ gives
\begin{equation}\label{eq:lambda_n}
\lambda_n=-\frac{1}{2}\,{\frac { \left( {q}^{n}-1 \right)  \left( 2\,\sqrt {q} \left({q}^{n}-q \right) a+ \left( q-1 \right)  \left( {q}^{n}+q \right) d
 \right) }{{q}^{n} \left( q-1 \right) ^{2}}}.
\end{equation}
Equating the coefficients of $B_n(\alpha,x)$ and $B_{n-1}(\alpha,x)$ gives $k'_n/k_n,~k''_n/k_n\in \mathbb{Q}(q^n,\sqrt{q})$.
\item Substitute the expression of $p_n$ given by \eqref{Pn} in the recurrence relation \eqref{TTRR} and use \eqref{ttrrxBN}.  By equating the coefficients of $B_{n+1}(\alpha,x)$, $B_{n}(\alpha,x)$, $B_{n-1}(\alpha,x)$, we get $A_n$, $B_n$ and $C_n$, respectively, given as
\begin{align}
  A_n\frac{ k_n}{ k_{n+1}}= &\frac{1}{\mu_2(\alpha,n)}, ~
  B_n\frac{ k_n}{ k_{n+1}}= -\frac{\mu_2(\alpha, n-1)k'_n}{\mu_2(\alpha,n)k_n}+\frac{k'_{n+1}}{k_{n+1}}-\frac{\mu_1(\alpha, n)}{\mu_2(\alpha, n)},\label{AnBnCn}\\
  C_n\frac{ k_{n-1}}{ k_{n+1}}=&-\frac{\mu_2(\alpha, n-1)(k'_n)^2}{\mu_2(\alpha, n)k_n^2}+\frac{\mu_2(\alpha, n-2)k''_n}{\mu_2(\alpha, n)k_n}-\frac{k''_{n+1}}{k_{n+1}}+\frac{k'_nk'_{n+1}}{k_nk_{n+1}}-\frac{(\mu_1(\alpha,n)-\mu_1(\alpha, n-1))k'_n}{\mu_2(\alpha, n)k_n}.\nonumber
\end{align}
  \item Substituting the values of $k'_n$ and $k''_n$ given in step 3 in these equations yields the three unknowns in terms of $\alpha$, $a$, $b$, $c$, $d$, $e$, $n$, $k_{n-1}$, $ k_n$, $ k_{n+1}$ given by $(N=q^n)$:
\end{enumerate}
\begin{align*}
&\frac{k_n}{k_{n+1}}A_n=-2\alpha N,\\
   &\frac{k_n}{k_{n+1}}B_n=2\,{N}^{2}\alpha\,\Big\{\left( -4\,{q}^{34} \left( q+1 \right)  \left( q-1 \right) ^{6}
 \left( N-1 \right)  \left( N-q \right) b+2\,{q}^{{\frac{67}{2}}}
 \left( q-1 \right) ^{8} \left( N+1 \right)  \left( N+q \right) e
 \right) a\\
 &-2\,{q}^{{\frac{67}{2}}} \left( q+1 \right)  \left( q-1
 \right) ^{7} \left( N-1 \right)  \left( N+q \right) db+{q}^{33}
 \left( q-1 \right) ^{8} \left( {N}^{2}q-N{q}^{2}-{N}^{2}-2\,Nq-{q}^{2
}-N+q \right) ed   \Big\}  \\
   &\Big/ \Big\{ 4\,{q}^{34}\left( q-1 \right) ^{6} \left( {N}^{2}-1 \right)
 \left( {N}^{2}-{q}^{2} \right) {a}^{2}+4\,{q}^{{\frac{67}{2}}}
 \left( q-1 \right)^7  \left( {N}^{4}-{q}^{2} \right) da+{q}^{33} \left( q-1 \right) ^{8} \left( {N}^{2}+1 \right)  \left( {N}^{2}+{q}^
{2} \right) {d}^{2}   \Big\}.
\end{align*}
The expression of $\frac{k_{n-1}}{k_{n+1}}C_n$ is very huge and will not be displayed here. However, it can be found in the Maple file associated to this manuscript at \url{http://www.mathematik.uni-kassel.de/~tcheutia/}.
\subsection{Polynomials expanded in the basis $\{\xi_n(\gamma,\delta,\mu(x))\}$}\label{TTRR_q_Racah}
The polynomials represented in this basis are:\\
1. the $q$--Racah polynomials
\[R_n(\mu(x);\alpha,\beta,\gamma,\delta|q)=\qhypergeom{4}{3}{q^{-n},\alpha\beta q^{n+1},q^{-x},\gamma\delta q^{x+1}}{\alpha q,\beta\delta q,\gamma q}{q;q},~n=0,1,\ldots,N,\]
with
\begin{eqnarray*}
  \phi(\mu(x)) &=& \frac { \left( q-1 \right) ^{2}}{2{q}^{3/2}}\Big( \left( {q}^{2}\alpha\,\beta+1
 \right)\mu(x)^2 - \left( q\gamma\,\delta\,\beta+q\delta
\,\beta\,\alpha+q\gamma\,\alpha+q\beta\,\alpha+\gamma+\gamma\,\delta+
\delta\,\beta+\alpha \right) \mu(x) \\
   && +2\left( q\gamma\,\delta\,\alpha+q\delta\,\beta
\,\alpha+q\gamma\,\alpha+q{\delta}^{2}\beta\,\gamma+q{\gamma}^{2}
\delta+q\gamma\,\delta\,\beta-\gamma\,\delta\,{q}^{2}\alpha\,\beta-
\gamma\,\delta \right)\Big), \\
  \psi(\mu(x)) &=& \frac { \left( q-1 \right)}{q}\Big(  \left( {q}^{2}\alpha\,\beta-1 \right) \mu(x) - q  \left( q\gamma\,\delta\,
\beta+q\delta\,\beta\,\alpha+q\gamma\,\alpha+q\beta\,\alpha-\gamma-
\gamma\,\delta-\delta\,\beta-\alpha \right)\Big);
\end{eqnarray*}
2. the dual $q$--Hahn polynomials
\[R_n(\mu(x);\gamma,\delta,N|q)=\qhypergeom{3}{2}{q^{-n},q^{-x},\gamma\delta q^{x+1}}{\gamma q,q^{-N}}{q;q},~n=0,1,\ldots,N,\]
with
\begin{eqnarray*}
  \phi(\mu(x)) &=& \frac { \left( q-1 \right) ^{2}}{2{q}^{3/2}}\Big(  \mu(x)^2- \frac{ \gamma\,q+{q}^{N}\gamma\,q+\gamma\,\delta\,{q}^{N}q+1 }{q^N}
\mu(x)+2\frac {\gamma q\, \left( \gamma\,\delta\,{q}^{N}q+1+\delta-{q}^{N}\delta
 \right) }{{q}^{N}}\Big), \\
  \psi(\mu(x))&=& {\frac { \left( q-1 \right)}{q}}\Big(- \mu(x) +\frac{ -\gamma\,q+{q}^{N}\gamma\,q+\gamma\,\delta\,{q}^{N
}q+1  }{{q}^{N}}\Big) ;
\end{eqnarray*}
3. the  dual $q$--Krawtchouk polynomials
\[K_n(\mu(x);c,N|q)=\qhypergeom{3}{2}{q^{-n},q^{-x},cq^{x-N}}{q^{-N},0}
{q;q},~n=0,1,\ldots,N,\]
with
\begin{eqnarray*}
  \phi(\mu(x)) &=& \frac{(q-1)^2}{2q^{3\over 2}}\Big(\mu(x)^2-\frac{c+1}{q^N}\mu(x)+2c\frac{1-q^N}{q^{2N}} \Big),\ \psi(\mu(x)) =\frac{q-1}{q}\Big(-\mu(x)+\frac{c+1}{q^N} \Big);
\end{eqnarray*}
 By direct computations, we have  the structure relations
\begin{align*}
  \xi_1(\gamma,\delta,\mu(x))\dx^2\xi_n(\gamma,\delta,\mu(x)) &= \eta(1,n)\eta(\sqrt{q},n-1)\xi_{n-1}(\gamma,\delta,\mu(x)); \\
  \xi_1(\gamma,\delta,\mu(x))\sx\dx \xi_n(\gamma,\delta,\mu(x)) &= \eta(1,n)\left(\beta_1(\sqrt{q},\gamma,\delta,n-1)\xi_{n-1}(\gamma,\delta,\mu(x))+\beta_2(n-1)\xi_n(\gamma,\delta,\mu(x))\right); \\
  \xi_1(\gamma,\delta,\mu(x))\xi_n(\gamma,\delta,\mu(x)) &= \nu_1(\gamma,\delta,n)\xi_n(\gamma,\delta,\mu(x))+\nu_2(n)\xi_{n+1}(\gamma,\delta,\mu(x));\\
  \mu(x)\xi_n(\gamma,\delta,\mu(x)) &= \mu_1(\gamma,\delta,n)\xi_n(\gamma,\delta,\mu(x))+\mu_2(n)\xi_{n+1}(\gamma,\delta,\mu(x));
\end{align*}
where
\begin{align*}
  &\mu_1(\gamma,\delta,n) = \frac{1+\gamma\delta q^{2n+1}}{q^n},\ \ \mu_2(n)=\frac{-1}{q^n},\  \nu_1(\gamma,\delta,n) = (1-q^{-n})(1-\gamma\delta q^{n+1}),\\
  &\nu_2(n)=q^{-n},~\beta_1(a,\gamma,\delta,n) = \frac{1}{2}(1-a^2\gamma\delta q^{2n})(1-q^{-n}),\ \ \beta_2(n)=\frac{1+q^n}{2q^n}, \ \ \eta(a,n)=\frac{a(1-q^n)}{q-1}.
\end{align*}
We suppose now that
\[ p_n(x):=p_n(\mu(x))=k_n\xi_n(\gamma,\delta,\mu(x))+k'_n\xi_{n-1}(\gamma,\delta,\mu(x))+k''_{n}\xi_{n-2}(\gamma,\delta,\mu(x))+\ldots.\]
To get the coefficients $A_n$, $B_n$ and $C_n$ of \eqref{TTRR}, we proceed as in section \ref{sect_Askey_Wilson} and obtain $\lambda_n$ given by \eqref{eq:lambda_n} and for $N=q^n$,
\begin{align*}
&\frac{k_n}{k_{n+1}}A_n=-N,\\
   &\frac{k_n}{k_{n+1}}B_n=\Big\{\left( -4\,{N}^{2}{q}^{4} \left( q+1 \right)  \left( q-1 \right) ^{2}
 \left( N-1 \right)  \left( N-q \right) b+2\,{N}^{2}{q}^{7/2} \left( q
-1 \right) ^{4} \left( N+1 \right)  \left( N+q \right) e \right) a\\
&-2\,{N}^{2}{q}^{7/2} \left( q+1 \right)  \left( q-1 \right) ^{3} \left( N-
1 \right)  \left( N+q \right) db+{N}^{2}{q}^{3} \left( q-1 \right) ^{4
} \left( {N}^{2}q-N{q}^{2}-{N}^{2}-2\,Nq-{q}^{2}-N+q \right) ed
 \Big\}   \\
   &\Big/ \Big\{  4\,{q}^{4} \left( {N}^{2}-1 \right)
 \left( {N}^{2}-{q}^{2} \right)\left( q-1 \right) ^{2} {a}^{2}+4\,{q}^{7/2} \left( q-1
 \right)^3  \left( {N}^{4}-{q}^{2} \right) da+{q}^{3} \left( q-1
 \right) ^{4} \left( {N}^{2}+1 \right)  \left( {N}^{2}+{q}^{2}
 \right) {d}^{2}
 \Big\}.
\end{align*}
As in  section \ref{sect_Askey_Wilson}, the coefficients $\frac{k_{n-1}}{k_{n+1}}C_n$ is huge and  can be found in the Maple file associated to this manuscript at \url{http://www.mathematik.uni-kassel.de/~tcheutia/}.
\section{Extension of the algorithms implemented in the Maple package \texttt{retode}}
As shown in section \ref{sect:TTRR}, the classical orthogonal polynomials on a quadratic or a $q$-quadratic lattice satisfy a recurrence equation
\[p_n(x)=(A_nx+B_n)p_n(x)-C_np_{n-1}(x),\]
with $A_n$, $B_n$, $C_n$ given explicitly. Koepf and Schmersau \cite{Koepf_Schmersau_2002}  implemented algorithms to test wether or not a given holonomic recurrence equation (i.~e. linear, homogeneous with polynomial coefficients)
\begin{equation}\label{ttrrg}
  q_n(x)p_{n+2}(x)+r_n(x)p_{n+1}(x)+s_n(x)p_n(x)=0~(q_n(x), \ r_n(x),\ s_n(x) \in\mathbb{Q}[n,x]),
\end{equation}
has  classical orthogonal polynomial solutions  of a continuous, a discrete or a $q$-discrete variable. In fact, it is less obvious to find out whether there is a polynomial system satisfying \eqref{ttrrg}, being a linear transform of one of the classical system (of a continuous, a discrete or a $q$-discrete variable), and to identify the system in the affirmative case.   In this section, we aim to implement, using the same approach, algorithms to test wether a given holonomic recurrence equation  has classical orthogonal polynomial solutions of a quadratic or a $q$-quadratic lattice. The algorithms were  explicitly given and explained in \cite{Koepf_Schmersau_2002} for classical orthogonal polynomials of a continuous, a discrete or a $q$-discrete variable and we will adapt them here for classical orthogonal polynomials on a quadratic or a $q$-quadratic lattice  according to the basis in which the polynomials are expanded.
\subsection{Polynomials expanded in the basis $\{\vartheta_n(\alpha,x)\}$}\label{algo_Wilson}
{\bf Algorithm 1} (see \cite[Algorithms 1 and 2]{Koepf_Schmersau_2002}). This algorithm decides whether a given holonomic three-term recurrence equation has classical orthogonal polynomial solutions  expanded in the basis $\{\vartheta_n(\alpha,x)\}$, and returns their data if applicable.
\begin{enumerate}
  \item \texttt{Input}:  A holonomic three-term recurrence equation
  \[q_n(x)p_{n+2}(x)+r_n(x)p_{n+1}(x)+s_n(x)p_n(x)=0~(q_n(x), \ r_n(x),\ s_n(x) \in\mathbb{Q}[n,x]).\]
  \item \texttt{Shift}: Shift by $\max\{n\in\mathbb{N}_0\ | \ n \text{ is zero of either } q_{n-1}(x) \text{ or } s_n(x) \}+1$ if necessary.
  \item \texttt{Rewriting}: Rewrite the recurrence equation in the form
  \[p_{n+1}(x)=t_n(x)p_n(x)+u_n(x)p_{n-1}(x)~~ (t_n(x),u_n(x)\in \mathbb{Q}(n,x)).\]
  If either $t_n(x)$ is not a polynomial of degree one in $x$ or $u_n(x)$ is not a constant with respect to $x$, return ``\texttt{no classical orthogonal polynomial solution exists}"; exit.
  \item \texttt{Linear transformation}: Rewrite the recurrence equation by the linear transformation $x\mapsto (x-g)/f$ with unknowns $f$ and $g$.
  \item\label{standardization} \texttt{Standardization}: Given now $A_n$, $B_n$ and $C_n$ by
  \[p_{n+1}(x)=(A_nx+B_n)p_n(x)-C_np_{n-1}(x)~(A_n,\ B_n,\ C_n\ \in\mathbb{Q}(n),\ A_n\neq 0),\]
  define
  \[\frac{k_{n+1}}{k_n}:=A_n=\frac{v_n}{w_n}~(v_n,\ w_n\ \in \mathbb{Q}[n])\]
  according to proposition \ref{prop1}.
  \item\label{makemonic} \texttt{Make monic}: Set
  \[\tilde{B}_n:=\frac{B_n}{A_n}\in\mathbb{Q}(n)\ \text{ and }\ \tilde{C}_n:=\frac{C_n}{A_nA_{n-1}}\in\mathbb{Q}(n)\]
  and bring these rational functions in lowest terms. If the degree of either the numerator $\tilde{B}_n$ is larger than 4, if the degree of  the denominator of $\tilde{B}_n$ is larger than 2, if the degree of the numerator of $\tilde{C}_n$ is larger than 8, or if the degree of the denominator of $\tilde{C}_n$ is larger than 4, then return ``\texttt{no classical orthogonal polynomial solution exists}".
  \item\label{poly_identities} \texttt{Polynomial identities}: Set
  \[\tilde{B}_n=\frac{k_n}{k_{n+1}}B_n,~ \tilde{C}_n=\frac{k_{n-1}}{k_{n+1}}C_n\]
  according to proposition \ref{prop1}, in terms of the unknowns $a,\ b,\ c,\ d,\ e,\ f$ and $g$. Multiply these identities by their common denominators, and bring them therefore  in polynomial form.
  \item \texttt{Equating coefficients}:  Equate the coefficients of the powers of $n$ in the two resulting equations. This results in a nonlinear system in the unknowns $a,\ b,\ c,\ d,\ e,\ f$ and $g$. Solve this system by Gr\"{o}bner bases methods. If the system has no solution or only one with $a=d=0$, then return `\texttt{no classical orthogonal polynomial solution exists}"; exit.
  \item \texttt{Output}:  Return the solution vectors $(a,b,c,d,f,g)$ of the last step, the  divided-difference equation \eqref{eqWCDH}  together with the information $\frac{k_{n+1}}{k_n}$ and $y=fx+g$.
\end{enumerate}

\begin{example}
For the first example, we consider the three-term recurrence equation satisfied by the Wilson polynomials
\[\tilde{W}_n(x^2):=\tilde{W}_n(x^2;a,b,c,d)=\frac{W_n(x^2;a,b,c,d)}{(a+b)_n(a+c)_n(a+d)_n} \]
 given by \cite[eq. (9.1.4)]{KLS}
\[-(a^2+x^2)\tilde{W}_n(x^2)=A_n\tilde{W}_{n+1}(x^2)-(A_n+C_n)\tilde{W}_n(x^2)+C_n\tilde{W}_{n-1}(x^2),\]
where
\begin{align*}
  A_n= & \frac { \left( n+a+b+c+d-1 \right)  \left( n+a+b \right)  \left( n+a+c \right)  \left( n+a+d \right) }{ \left( 2\,n+a+b+c+d-1 \right) \left( 2\,n+a+b+c+d \right) }, \\
  C_n= & \frac {n \left( n+b+c-1 \right)  \left( n+b+d-1 \right)  \left( n+c+d-1 \right) }{ \left( 2\,n+a+b+c+d-2 \right)  \left( 2\,n+a+b+c+d-1
 \right) }.
\end{align*}
Using our implementation, the result is obtained by
\begin{maplegroup}
\begin{mapleinput}
\mapleinline{active}{1d}{A[n]:=(n+a+b+c+d-1)*(n+a+b)*(n+a+c)*(n+a+d)/((2*n+a+b+c+d-1)*(2*n+a+b+c+d))
}{}
\end{mapleinput}
\mapleresult
\begin{maplelatex}
\mapleinline{inert}{2d}{A[n] := (n+a+b+c+d-1)*(n+a+b)*(n+a+c)*(n+a+d)/((2*n+a+b+c+d-1)*(2*n+a+b+c+d))}{\[\displaystyle A_{{n}}\, := \,{\frac { \left( n+a+b+c+d-1 \right)  \left( n+a+b \right)  \left( n+a+c \right)  \left( n+a+d \right) }{ \left( 2\,n+a+b+c+d-1 \right) \\
\mbox{} \left( 2\,n+a+b+c+d \right) }}\]}
\end{maplelatex}
\end{maplegroup}
\begin{maplegroup}
\begin{mapleinput}
\mapleinline{active}{1d}{C[n]:=n*(n+b+c-1)*(n+b+d-1)*(n+c+d-1)/((2*n+a+b+c+d-2)*(2*n+a+b+c+d-1))
}{}
\end{mapleinput}
\mapleresult
\begin{maplelatex}
\mapleinline{inert}{2d}{C[n] := n*(n+b+c-1)*(n+b+d-1)*(n+c+d-1)/((2*n+a+b+c+d-2)*(2*n+a+b+c+d-1))}{\[\displaystyle C_{{n}}\, := \,{\frac {n \left( n+b+c-1 \right)  \left( n+b+d-1 \right)  \left( n+c+d-1 \right) }{ \left( 2\,n+a+b+c+d-2 \right)  \left( 2\,n+a+b+c+d-1 \right) \\
\mbox{}}}\]}
\end{maplelatex}
\end{maplegroup}
\begin{maplegroup}
\begin{mapleinput}
\mapleinline{active}{1d}{RecWilson:= -(a\symbol{94}2+x)*p(n)=A[n]*p(n+1)-(A[n]+C[n])*p(n)+C[n]*p(n-1):
}{}
\end{mapleinput}
\end{maplegroup}
\begin{maplegroup}
\begin{mapleinput}
\mapleinline{active}{1d}{strict:=true:
}{}
\end{mapleinput}
\end{maplegroup}
\begin{maplegroup}
\begin{mapleinput}
\mapleinline{active}{1d}{REtoWilsonDE(subs(n=n+1, RecWilson),p(n),x)
}{}
\end{mapleinput}
\mapleresult
\begin{maplelatex}
\mapleinline{inert}{2d}{`Warning: parameters have the values`, {a = a, b = -a*(a*b+a*c+a*d+b*c+b*d+c*d),}}{\[\displaystyle \mbox {{\tt `Warning: parameters have the values`}},\, \Big\{ a=a,b=-a \left( ab+ac+ad+cb+bd+cd \right) , \]}
\end{maplelatex}
\begin{maplelatex}
\mapleinline{inert}{2d}{ { c = (a*a)*b*c*d, d = a*a+a*b+a*c+a*d, e = -(a*a)*b*c-(a*a)*b*d-(a*a)*c*d-a*b*c*d}}{\[\displaystyle    c={a}^{2}bcd,d={a}^{2}+ab+ac+ad,e=-{a}^{2}bc-{a}^{2}bd-{a}^{2}cd-abcd \Big\} \]}
\end{maplelatex}
\mapleresult
\begin{maplelatex}
\mapleinline{inert}{2d}{[(a*b*c*d-a*b*x-a*c*x-a*d*x-b*c*x-b*d*x-c*d*x+x^2)*DD(DD(p(n, x), x), x)-(a*b*c+a*b*d+a*c*d+b*c*d-a*x-b*x-c*x-d*x)*SS(DD(p(n, x), x), x)-n*(n+a+b+c+d-1)*p(n, x) = 0, k[n+1]/k[n] = -(2*n+a+b+c+d)*(2*n+a+b+c+d-1)/((n+a+d)*(n+a+c)*(n+a+b)*(n+a+b+c+d-1))]}{\[\displaystyle \Big[ \left( abcd-abx-acx-adx-cbx-bdx-cdx+{x}^{2} \right) {\it DD} \left( {\it DD} \left( p \left( n,x \right) ,x \right) ,x \right) \]}
\end{maplelatex}
\begin{maplelatex}
\mapleinline{inert}{2d}{[(a*b*c*d-a*b*x-a*c*x-a*d*x-b*c*x-b*d*x-c*d*x+x^2)*DD(DD(p(n, x), x), x)-(a*b*c+a*b*d+a*c*d+b*c*d-a*x-b*x-c*x-d*x)*SS(DD(p(n, x), x), x)-n*(n+a+b+c+d-1)*p(n, x) = 0, k[n+1]/k[n] = -(2*n+a+b+c+d)*(2*n+a+b+c+d-1)/((n+a+d)*(n+a+c)*(n+a+b)*(n+a+b+c+d-1))]}{\[\displaystyle - \left( abc+abd+acd+bcd-ax-bx-cx-dx \right) {\it SS} \left( {\it DD} \left( p \left( n,x \right) ,x \right) ,x \right)\]}
\end{maplelatex}
\begin{maplelatex}
\mapleinline{inert}{2d}{[(a*b*c*d-a*b*x-a*c*x-a*d*x-b*c*x-b*d*x-c*d*x+x^2)*DD(DD(p(n, x), x), x)-(a*b*c+a*b*d+a*c*d+b*c*d-a*x-b*x-c*x-d*x)*SS(DD(p(n, x), x), x)-n*(n+a+b+c+d-1)*p(n, x) = 0, k[n+1]/k[n] = -(2*n+a+b+c+d)*(2*n+a+b+c+d-1)/((n+a+d)*(n+a+c)*(n+a+b)*(n+a+b+c+d-1))]}{\[\displaystyle-n \left( n+a+b+c+d-1 \right) p \left( n,x \right) =0,\]}
\end{maplelatex}
\begin{maplelatex}
\mapleinline{inert}{2d}{[(a*b*c*d-a*b*x-a*c*x-a*d*x-b*c*x-b*d*x-c*d*x+x^2)*DD(DD(p(n, x), x), x)-(a*b*c+a*b*d+a*c*d+b*c*d-a*x-b*x-c*x-d*x)*SS(DD(p(n, x), x), x)-n*(n+a+b+c+d-1)*p(n, x) = 0, k[n+1]/k[n] = -(2*n+a+b+c+d)*(2*n+a+b+c+d-1)/((n+a+d)*(n+a+c)*(n+a+b)*(n+a+b+c+d-1))]}{\[\displaystyle {\frac {k_{{n+1}}}{k_{{n}}}}=-{\frac { \left( 2\,n+a+b+c+d-1 \right)  \left( 2\,n+a+b+c+d \right) \\
\mbox{}}{ \left( n+a+b+c+d-1 \right)  \left( n+a+b \right)  \left( n+a+c \right)  \left( n+a+d \right) }}\Big]\]}
\end{maplelatex}
\end{maplegroup}
\noindent which gives the divided-difference equation of the Wilson polynomials (see \cite[theo. 2.5]{njionou_et_al_2015}), as well as the term ratio $k_{n+1}/k_n$. Here \text{SS} and \text{DD} stand for $\mathbf{S}_x$ and $\mathbf{D}_x$, respectively.
\end{example}
\begin{example}
  Abdulaziz D. Alhaidari \cite{Alhaidari} encountered two families of orthogonal polynomials on the real line defined by their three-term recurrence relations and initial values. The first system is given by
  \begin{align}
    &\cos\theta H_n^{(\mu,\nu)}(z;\alpha, \theta)=\left(z\sin\theta \left[\left(n+\frac{\mu+\nu+1}{2} \right)^2+\alpha \right]+\frac{\nu^2-\mu^2}{(2n+\mu+\nu)(2n+\mu+\nu+2)}  \right)H_n^{(\mu,\nu)}(z;\alpha, \theta) \nonumber \\
     &+\frac{2(n+\mu)(n+\nu)}{(2n+\mu+\nu)(2n+\mu+\nu+1)}H_{n-1}^{(\mu,\nu)}(z;\alpha, \theta)+\frac{2(n+1)(n+\mu+\nu+1)}{(2n+\mu+\nu+1)(2n+\mu+\nu+2)}H_{n+1}^{(\mu,\nu)}(z;\alpha, \theta),\label{Alhaidari1}
  \end{align}
  with $0\leq \theta \leq \pi$, $\mu,\nu>-1$, $\alpha\in\mathbb{R}$ and initial values $H_0^{(\mu,\nu)}(z;\alpha, \theta)=1$, $H_{-1}^{(\mu,\nu)}(z;\alpha, \theta)=0$. \\
  The second system is
  \begin{align}
   zG_n^{(\mu, \nu)}(z;\sigma)&=\left((\sigma+B_n^2)\left[\frac{\mu^2-\nu^2}{(2n+\mu+\nu)(2n+\mu+\nu+2)}+1 \right]-\frac{2n(n+\nu)}{2n+\mu+\nu}-\frac{(\mu+1)^2}{2} \right)G_n^{(\mu,\nu)}(z;\sigma)\nonumber \\
     & -(\sigma+B_{n-1}^2)\frac{2(n+\mu)(n+\nu)}{(2n+\mu+\nu)(2n+\mu+\nu+1)}G_{n-1}^{(\mu,\nu)}(z;\sigma)\nonumber\\
     &-(\sigma+B_n^2)\frac{2(n+1)(n+\mu+\nu+1)}{(2n+\mu+\nu+1)(2n+\mu+\nu+2)}G_{n+1}^{(\mu,\nu)}(z;\sigma),\label{Alhaidari2}
  \end{align}
  with $B_n=n+1+\frac{\mu+\nu}{2}$, $\mu,\nu>-1$ and $\sigma\in \mathbb{R}$ and initial values $G_0^{(\mu,\nu)}(z;\sigma)=1$, $G_{-1}^{(\mu,\nu)}(z;\sigma)=0$.
\end{example}
\noindent  Our implementation with
\begin{maplegroup}
\begin{mapleinput}
\mapleinline{active}{1d}{BB:=n->n+1+(mu+nu)/2:
}{}
\end{mapleinput}
\end{maplegroup}
\begin{maplegroup}
\begin{mapleinput}
\mapleinline{active}{1d}{recOpen2:=z*p(n)=((sigma+BB(n)\symbol{94}2)*((mu\symbol{94}2-nu\symbol{94}2)/((2*n+mu+nu)*(2*n+mu+nu+2))+1)}{}
\end{mapleinput}
\begin{mapleinput}
\mapleinline{active}{1d}{-2*n*(n+nu)/(2*n+mu+nu)-(mu+1)\symbol{94}2/2)*p(n)-(sigma+BB(n-1)\symbol{94}2)*2*(n+mu)*(n+nu)
/((2*n+mu+nu)*(2*n+mu+nu+1))*p(n-1)
}{}
\end{mapleinput}
\begin{mapleinput}
\mapleinline{active}{1d}{-(sigma+BB(n)\symbol{94}2)*2*(n+1)*(n+mu+nu+1)/((2*n+mu+nu+1)*(2*n+mu+nu+2))*p(n+1):
}{}
\end{mapleinput}
\end{maplegroup}
\begin{maplegroup}
\begin{mapleinput}
\mapleinline{active}{1d}{strict:=false:
}{}
\end{mapleinput}
\end{maplegroup}
\begin{maplegroup}
\begin{mapleinput}
\mapleinline{active}{1d}{REtoWilsonDE(subs(n=n+1, recOpen2),p(n),z)
}{}
\end{mapleinput}
\end{maplegroup}
\noindent returns six divided-difference equations (see the  Maple file associated to this manuscript) for the second recurrence equation \eqref{Alhaidari2}:
\begin{align}
   &\left({z}^{2}+ \left( {\mu}^{2}-2\,\sigma-3 \right) z+\frac{1}{4}\, \left( \mu-1 \right) ^{2} \left( {\mu}^{2}+2\,\mu+4\,\sigma+1 \right)
 \right)\mathbf{D}_z^2P_n(z/2) \nonumber \\
   &+(4\,z+2\, \left( \mu-1 \right)  \left( \mu+2\,\sigma+1 \right))\mathbf{S}_z\mathbf{D}_z P_n(z/2)-4n(n+1)P_n(z/2)=0,\label{Ald1}
\end{align}
\begin{align}
   &\left({z}^{2}+ \left({\mu}^{2}+2\,\mu-2\,\sigma+1 \right) z+\frac{1}{4} \left( \mu+1 \right) ^{2} \left( {\mu}^{2}+2\,\mu+4\,\sigma+1
 \right) \right)\mathbf{D}_z^2P_n(z/2) \nonumber \\
  & -4\,\sigma\, \left( \mu+1 \right)\mathbf{S}_z\mathbf{D}_z P_n(z/2)-4n(n-1)P_n(z/2)=0,\label{Ald2}
\end{align}
\begin{align}
   & \Big({a}^{2}{z}^{2}+ \left(  \left( {\mu}^{2}+2\,\mu-2\,\sigma+1 \right) {a
}^{2}-d \left( \mu+1 \right) a-\frac{1}{2}\,{d}^{2} \right) z\nonumber \\
&+\frac{1}{4}\, \left( \mu+1 \right) ^{2} \left(  \left( {\mu}^{2}+2\,\mu+4\,\sigma+1 \right) {a}^{2}-2\,d \left( \mu+1 \right) a+{d}^{2} \right)  \Big)\mathbf{D}_z^2P_n(z/2) \nonumber \\
   &+\left(2\,adz- \left( \mu+1 \right)  \left( 4\,{a}^{2}\sigma-d \left( \mu+1 \right) a+{d}^{2} \right) \right)\mathbf{S}_z\mathbf{D}_z P_n(z/2)\nonumber\\
   &-4na((n-1)a+d)P_n(z/2)=0~(a\neq 0,~ d\neq 0),\label{Ald3}
\end{align}
\begin{align}
   & (\left( 4\,\nu+1 \right) z-\frac{1}{2}\,{\nu}^{2})\mathbf{D}_z^2P_n(2z)+(\nu^2+\nu-2z)\mathbf{S}_z\mathbf{D}_z P_n(2z)+nP_n(2z)=0~(\nu\neq -\frac{1}{4}),\label{Ald4}
\end{align}
with $\sigma=0$, $\mu=\nu-1$,
\begin{align}
   &\left(\left( 4\,\nu+3 \right) z-\frac{1}{2}\,\nu\, \left( \nu+1 \right) \right)\mathbf{D}_z^2P_n(2z)+(\nu^2+2\nu-2z+\frac{1}{2})\mathbf{S}_z\mathbf{D}_z P_n(2z)+nP_n(2z)=0~(\nu\neq -\frac{3}{4}),\label{Ald5}
\end{align}
with $\sigma=-\frac{1}{4}$, $\mu=\nu$,
\begin{align}
   &\left(\left( 4\,\nu+5 \right) z-\frac{1}{2}\, \left( \nu+1 \right)^{2} \right)\mathbf{D}_z^2P_n(2z)+(\nu^2+3\nu-2z+2)\mathbf{S}_z\mathbf{D}_z P_n(2z)+nP_n(2z)=0~(\nu\neq -\frac{5}{4}),\label{Ald6}
\end{align}
with $\sigma=0$, $\mu=\nu+1$. By comparison with the Wilson divided-difference equation, we deduce from  the first three divided-difference equations \eqref{Ald1}--\eqref{Ald3}  that
\[G_n^{(\mu, \nu)}(z;\sigma)=constant\times W_n(z/2;a,b,c,d)\]
  where
$a,b,c,d$  are permutations of elements of the set   $\{ \frac{1}{2}(-\mu+1), \frac{1}{2}(-\mu+1), \frac{1}{2}(\mu+1)+\sqrt{-\sigma},  \frac{1}{2}(\mu+1)-\sqrt{-\sigma}\}$, $\{ \frac{1}{2}(\mu+1), \frac{1}{2}(\mu+1), -\frac{1}{2}(\mu+1)+\sqrt{-\sigma},  -\frac{1}{2}(\mu+1)-\sqrt{-\sigma}\}$,  $\{ \frac{1}{2}(\mu+1), \frac{1}{2}(\mu+1), \frac{1}{2}(\delta-\mu-1)+\sqrt{-\sigma},  \frac{1}{2}(\delta-\mu-1)-\sqrt{-\sigma}\}$ for the first equation \eqref{Ald1}, the second equation \eqref{Ald2}, and
 the third equation \eqref{Ald3} in which the parameters $a=1$, $d=\delta$, respectively. This brings then a new parameter $\delta$ in the definition of the polynomial $G_n^{(\mu, \nu)}(z;\sigma)$ and we also remark that for $d=0$ and $a=1$ in the third equation \eqref{Ald3}, we recovered  the second equation \eqref{Ald2}. For the value $d=\delta=n+\mu+2$, we recover (from the identification of \eqref{Ald3} with the Wilson polynomials) the solution
\[G_n^{(\mu, \nu)}(z;\sigma)=\frac{ W_n(z/2;a,b,c,d)}{(a+b)_n(a+d)_n}\]  given in \cite{Clarkson_van_Assche} where $a=c=\frac{\mu+1}{2}$, $b=\frac{n+1}{2}+\sqrt{-\sigma}$ and $d=\frac{n+1}{2}-\sqrt{-\sigma}$.\\
Comparing  the last three equations \eqref{Ald4}--\eqref{Ald6} with the divided-difference equation of the continuous dual Hahn polynomials $S_n(2z;a,b,c)$, we deduce that
  \[G_n^{(\nu-1, \nu)}(z;0)=constant\times S_n(2z;a,b,c),\]
  with $ a=c=\nu,~ b=\frac{1}{2}$, or $a=b=\nu, c=\frac{1}{2}$ or $b=c=\nu, a=\frac{1}{2}$, i.~e.,  $a,b,c$ are permutations of elements of the set $\{\nu, \nu, \frac{1}{2}\}$;
  \[G_n^{(\nu, \nu)}(z;-1/4)=constant\times S_n(2z;a,b,c),\]
  where $a,b,c$  are permutations of elements of the set   $\{\nu+1,~\nu,~\frac{1}{2}\}$;
  \[G_n^{(\nu+1, \nu)}(z;0)=constant\times S_n(2z;a,b,c),\]
  where $a,b,c$ are permutations of the elements of the set $\{\nu+1, \nu+1,\frac{1}{2}\}$.
\subsection{Polynomials expanded in the basis $\{(\alpha+ix)_n\}$}
The steps of the algorithm in this case agree with those given in section \ref{algo_Wilson}.  In steps \ref{standardization} and \ref{poly_identities}, we use proposition \ref{prop:CHMP} whereas in step \ref{makemonic}, the algorithm will return ``\texttt{no classical orthogonal polynomial solution exists}" if the degree of either the numerator or   the denominator of $\tilde{B}_n$ is larger than 2, if the degree of the numerator of $\tilde{C}_n$ is larger than 7, or if the degree of the denominator of $\tilde{C}_n$ is larger than 5.
\begin{example}
As example here, starting from the three-term recurrence equations (\texttt{RE}) \cite[eq. (9.4.3)]{KLS} and \cite[eq. (9.7.3)]{KLS} satisfied by the  continuous Hahn and the Meixner-Pollaczek polynomials, respectively, and using our implementation  with \texttt{REtoContHahnDE(subs(n=n+1, RE), p(n), x)}, we recover the divided-difference equations of type \eqref{CHMP} satisfied by both  families (see \cite[prop. 4]{Tcheutia_et_al_2017}). In the output in this case, \text{SS} and \text{DD} stand for $\mathbf{S}_x$ and $\mathbf{\delta}_x$, respectively.
\end{example}

\subsection{Polynomials expanded in the basis $\{\ \chi_n(\gamma,\delta,\lambda(x))\}$}
We proceed as in the algorithm of section \ref{algo_Wilson}. Here, in steps \ref{standardization} and \ref{poly_identities}, we use proposition \ref{prop_Racah} whereas in step \ref{makemonic}, the algorithm will return ``\texttt{no classical orthogonal polynomial solution exists}" if the degree of either the numerator of $\tilde{B}_n$ is larger than 4,  the degree of the denominator of $\tilde{B}_n$ is larger than 2, if the degree of the numerator of $\tilde{C}_n$ is larger than 8, or if the degree of the denominator of $\tilde{C}_n$ is larger than 4.
\begin{example}
  If we call \texttt{RE} the three-term recurrence equation of the Racah or  the dual Hahn polynomials given, respectively, by \cite[eq. (9.2.3)]{KLS} and \cite[eq. (9.6.3)]{KLS}, then with our implementation \texttt{REtoRacahDE(subs(n=n+1, RE), p(n),x)}, we get the divided-difference equation satisfied by both families with $\phi$ and $\psi$ given as in section \ref{sect:racah_DH}.
\end{example}
\begin{remark}
 From our implementations of sections 3.1. 3.2 and 3.3, we get for the recurrence equation \eqref{Alhaidari1} the following:
  \begin{maplegroup}
\begin{mapleinput}
\mapleinline{active}{1d}{RE:=1/2*(y+y\symbol{94}(-1))*p(n)=(z/(2*I)*(y-y\symbol{94}(-1))*((n+(mu+nu+1)/2)\symbol{94}2+alpha)+ (nu\symbol{94}2-mu\symbol{94}2)/((2*n+mu+nu)*(2*n+mu+nu+2)))*p(n)
+2*(n+mu)*(n+nu)/((2*n+mu+nu)*(2*n+mu+nu+1))*p(n-1)
+2*(n+1)*(n+mu+nu+1)/((2*n+mu+nu+1)*(2*n+mu+nu+2))*p(n+1):
}{}
\end{mapleinput}
\end{maplegroup}
\begin{maplegroup}
\begin{mapleinput}
\mapleinline{active}{1d}{strict:=false:
}{}
\end{mapleinput}
\end{maplegroup}
\begin{maplegroup}
\begin{mapleinput}
\mapleinline{active}{1d}{REtoWilsonDE(subs(n=n+1, RE), p(n), z)
}{}
\end{mapleinput}
\mapleresult
\begin{maplelatex}
\mapleinline{inert}{2d}{`Warning: parameters have the values`, {a = a, alpha = alpha, b = b, c = c, d = d, }}{\[\displaystyle \mbox {{\tt `Warning: parameters have the values`}},\, \left\{ a=a,\alpha=\alpha,b=b,c=c,d=d,\right\} \]}
\end{maplelatex}
\begin{maplelatex}
\mapleinline{inert}{2d}{ { e = e, f = 0, g = g, mu = mu, nu = nu, y = -1}}{\[\displaystyle \, \left\{ e=e,f=0,g=g,\mu=\mu,\nu=\nu,y=1 \right\} \]}
\end{maplelatex}
\mapleresult
\begin{maplelatex}
\mapleinline{inert}{2d}{`Warning: parameters have the values`, {a = a, alpha = alpha, b = b, c = c, d = d, }}{\[\displaystyle \mbox {{\tt `Warning: parameters have the values`}},\, \left\{ a=a,\alpha=\alpha,b=b,c=c,d=d,\right\} \]}
\end{maplelatex}
\begin{maplelatex}
\mapleinline{inert}{2d}{ { e = e, f = 0, g = g, mu = mu, nu = nu, y = 1}}{\[\displaystyle \, \left\{ e=e,f=0,g=g,\mu=\mu,\nu=\nu,y=1 \right\} \]}
\end{maplelatex}
\mapleresult
\begin{maplelatex}
\mapleinline{inert}{2d}{[[DD(DD(p(n, g), z), z)+(d*g+e)*SS(DD(p(n, g), z), z)/(a*g^2+b*g+c)-n*(a*n-a+d)*p(n, g)/(a*g^2+b*g+c) = 0],}{\[\displaystyle [[{\it DD} \left( {\it DD} \left( p \left( n,g \right) ,z \right) ,z \right) +{\frac { \left( dg+e \right) {\it SS} \left( {\it DD} \left( p \left( n,g \right) ,z \right) ,z \right) }{a{g}^{2}+bg+c}}\\
\mbox{}-{\frac {n \left( an-a+d \right) p \left( n,g \right) }{a{g}^{2}+bg+c}}=0],]\]}
\end{maplelatex}
\end{maplegroup}
\noindent We get the same answer using
\begin{maplegroup}
\begin{mapleinput}
\mapleinline{active}{1d}{REtoRacahDE(subs(n=n+1, RE),p(n),z)}{}
\end{mapleinput}
\end{maplegroup}
\noindent and
\begin{maplegroup}
\begin{mapleinput}
\mapleinline{active}{1d}{REtoContHahnDE(subs(n=n+1, RE),p(n),z)
}{}
\end{mapleinput}
\end{maplegroup}
 \noindent  We deduce from our implementations (which return the  solution $H_n^{(\mu,\nu)}(z;\alpha, \theta)=0$) that the polynomial family $H_n^{(\mu,\nu)}(z;\alpha, \theta)$ satisfying the recurrence equation \eqref{Alhaidari1} is not related to a known classical orthogonal polynomial sequence on a quadratic lattice expanded in the basis $\{\vartheta_n(\alpha,x)\}$, $\{(\alpha+ix)_n\}$ or  $\{\ \chi_n(\gamma,\delta,\lambda(x))\}$. This recurrence equation may lead to a new family of orthogonal polynomial sequence.
\end{remark}

\subsection{Polynomials expanded in the basis $\{B_n(\alpha,x)\}$}\label{algo_AW}
{\bf Algorithm 2} (see \cite[Algorithm 3]{Koepf_Schmersau_2002}). This algorithm decides whether a given holonomic three-term recurrence equation has classical orthogonal polynomial solutions  expanded in the basis $\{B_n(\alpha,x)\}$, and returns their data if applicable.
\begin{enumerate}
  \item \texttt{Input}:  A holonomic three-term recurrence equation
  \[q_n(x)p_{n+2}(x)+r_n(x)p_{n+1}(x)+s_n(x)p_n(x)=0~\left(q_n(x), \ r_n(x),\ s_n(x) \in\mathbb{Q}[q^n,\sqrt{q},x]\right).\]
  \item \texttt{Shift}: Shift by $\max\{n\in\mathbb{N}_0\ | \ n \text{ is zero of either } q_{n-1}(x) \text{ or } s_n(x) \}+1$ if necessary.
  \item \texttt{Rewriting}: Rewrite the recurrence equation in the form
  \[p_{n+1}(x)=t_n(x)p_n(x)+u_n(x)p_{n-1}(x)~~ \left(t_n(x),u_n(x)\in \mathbb{Q}(q^n,\sqrt{q},x)\right).\]
  If either $t_n(x)$ is not a polynomial of degree one in $x$ or $u_n(x)$ is not a constant with respect to $x$, return ``\texttt{no classical orthogonal polynomial solution exists}"; exit.
  \item \texttt{Linear transformation}: Rewrite the recurrence equation by the linear transformation $x\mapsto (x-g)/f$ with unknowns $f$ and $g$.
  \item\label{standardization1} \texttt{Standardization}: Given now $A_n$, $B_n$ and $C_n$ by
  \[p_{n+1}(x)=(A_nx+B_n)p_n(x)-C_np_{n-1}(x)~\left(A_n,\ B_n,\ C_n\ \in\mathbb{Q}(q^n,q),\ A_n\neq 0\right),\]
  define
  \[\frac{k_{n+1}}{k_n}:=-\frac{1}{2\alpha q^n}A_n=\frac{v_n}{w_n}~\left(v_n,\ w_n\ \in \mathbb{Q}[q^n, q]\right).\]
  \item\label{makemonic1} \texttt{Make monic}: Set
  \[\tilde{B}_n:=\frac{B_n}{A_n}\in\mathbb{Q}(q^n,\sqrt{q})\ \text{ and }\ \tilde{C}_n:=\frac{C_n}{A_nA_{n-1}}\in\mathbb{Q}(q^n,\sqrt{q})\]
  and bring these rational functions in lowest terms. If the degree (w.~r.~t.\ $N:=q^n$) of either the numerator $\tilde{B}_n$ is larger than 3, if the degree of  the denominator of $\tilde{B}_n$ is larger than 4, if the degree of the numerator  or the denominator of $\tilde{C}_n$ is larger than 14, then return ``\texttt{no classical orthogonal polynomial solution exists}".
  \item\label{poly_identities1} \texttt{Polynomial identities}: Set
  \[\tilde{B}_n=\frac{k_n}{k_{n+1}}B_n,~ \tilde{C}_n=\frac{k_{n-1}}{k_{n+1}}C_n\]
   in terms of the unknowns $a,\ b,\ c,\ d,\ e,\ f$ and $g$. Multiply these identities by their common denominators, and bring them therefore  in polynomial form.
  \item \texttt{Equating coefficients}:  Equate the coefficients of the powers of $N=q^n$ in the two resulting equations. This results in a nonlinear system in the unknowns $a,\ b,\ c,\ d,\ e,\ f$ and $g$. Solve this system by Gr\"{o}bner bases methods. If the system has no solution or only one with $a=d=0$, then return `\texttt{no classical orthogonal polynomial solution exists}"; exit.
  \item \texttt{Output}:  Return the solution vectors $(a,b,c,d,f,g)$ of the last step, the  divided-difference equation \eqref{eq:dxsx}  together with the information $\frac{k_{n+1}}{k_n}$ and $y=fx+g$.
\end{enumerate}
\begin{example}
As illustrative  example, we use our implementation to find the divided-difference equation of type \eqref{eq:dxsx}  satisfied by the continuous Hermite polynomials.
\begin{maplegroup}
\begin{mapleinput}
\mapleinline{active}{1d}{recContinuousqHermite:=2*x*p(n)=p(n+1)+(1-q\symbol{94}n)*p(n-1)
}{}
\end{mapleinput}
\mapleresult
\begin{maplelatex}
\mapleinline{inert}{2d}{recContinuousqHermite := 2*x*p(n) = p(n+1)+(1-q^n)*p(n-1)}{\[\displaystyle {\it recContinuousqHermite}\, := \,2\,xp \left( n \right) =p \left( n+1 \right) + \left( 1-{q}^{n} \right) p \left( n-1 \right) \]}
\end{maplelatex}
\end{maplegroup}
\begin{maplegroup}
\begin{mapleinput}
\mapleinline{active}{1d}{strict:=true:
}{}
\end{mapleinput}
\end{maplegroup}
\begin{maplegroup}
\begin{mapleinput}
\mapleinline{active}{1d}{REtoAskeyWilsonDE(subs(n=n+1, recContinuousqHermite), p(n),q, x)
}{}
\end{mapleinput}
\mapleresult
\begin{maplelatex}
\mapleinline{inert}{2d}{[(1/2*(2*x^2-1))*DD(DD(p(n, x), x), x)-2*x*SS(DD(p(n, x), x), x)*sqrt(q)/(q-1)]}{\[\displaystyle [1/2\, \left( 2\,{x}^{2}-1 \right) {\it DD} \left( {\it DD} \left( p \left( n,x \right) ,x \right) ,x \right) -2\,{\frac {x{\it SS} \left( {\it DD} \left( p \left( n,x \right) ,x \right) ,x \right)  \sqrt{q}}{q-1}}\]}
\end{maplelatex}
\begin{maplelatex}
\mapleinline{inert}{2d}{[+2*q^(3/2)*(-1+q^n)*p(n, x)/(q^n*(q-1)^2) = 0, k[n+1]/k[n] = 2]}{\[\displaystyle
\mbox{}+2\,{\frac {{q}^{3/2} \left( -1+{q}^{n} \right) p \left( n,x \right) }{{q}^{n} \left( q-1 \right) ^{2}}}=0,{\frac {k_{{n+1}}}{k_{{n}}}}=2]\]}
\end{maplelatex}
\end{maplegroup}
In the result, \text{SS} and \text{DD} stand for $\mathbb{S}_x$ and $\mathbb{D}_x$, respectively. The results for the other families can be found in the accompanying Maple file of this manuscript.
\end{example}

\subsection{Polynomials expanded in the basis $\{\xi_n(\gamma,\delta,\mu(x))\}$}
The steps of the algorithm in this case agree with those given in section \ref{algo_AW}.  In steps \ref{standardization1} and \ref{poly_identities1}, we use the results from section \ref{TTRR_q_Racah}  whereas in step \ref{makemonic1}, the algorithm will return ``\texttt{no classical orthogonal polynomial solution exists}" if the degree of either the numerator of $\tilde{B}_n$ is larger than 3, if   the denominator of $\tilde{B}_n$ is larger than 4, if the degree of the numerator  or  the denominator of $\tilde{C}_n$ is larger than 12.
\begin{example}
If we consider for example the recurrence equation \texttt{RE} for the $q$-Racah, dual $q$-Hahn, dual $q$-Krawtchouk given, respectively,  by \cite[eq. (14.2.3)]{KLS}, \cite[eq. (14.7.3)]{KLS}, \cite[eq. (14.17.3)]{KLS}, we use our implementation \texttt{REtoqRacahDE(subs(n=n+1,RE), p(n), x, q)} to get the divided-difference equations satisfied by the three families of polynomials and the product $\gamma\delta$.
\end{example}
\noindent \texttt{Note}: The Maple implementation \texttt{retode} by Koepf and Schmersau has been updated with our extension to classical orthogonal polynomials on a quadratic or a $q$-quadratic lattice. The package \texttt{retode.mpl} and a worksheet \texttt{retodedemo.mw} containing the three-term recurrence equations of section 2 and the examples for all the classical orthogonal polynomials on a quadratic or a $q$-quadratic lattice can be obtained from \url{http://www.mathematik.uni-kassel.de/~tcheutia/}.

\section{Acknowledgments}
This work has been supported by the Institute of Mathematics of the University of Kassel, Germany. The author thanks Wolfram Koepf for bringing his attention to this problem and for helpful discussions to improve his Maple implementation.

\end{document}